\documentclass[11pt,amsfonts,amssymb]{article}
\usepackage{amsfonts}
\usepackage{amssymb}

\newcommand{\qed}{\ \hfill \rule{1ex}{1ex}}

\newcommand{\ci}{\subseteq}

\newcommand{\union}{\cup}

\newcommand{\infinity}{\infty}

\newcommand{\bprop}{\begin{proposition}}
\newcommand{\eprop}{\end{proposition}}
\newcommand{\bthm}{\begin{theorem}}
\newcommand{\ethm}{\end{theorem}}
\newcommand{\blem}{\begin{lemma}}
\newcommand{\elem}{\end{lemma}}
\newcommand{\bpf}{\begin{proof}}
\newcommand{\epf}{\end{proof}}
\newcommand{\bdefn}{\begin{definition}}
\newcommand{\edefn}{\end{definition}}
\newcommand{\benum}{\begin{enumerate}}
\newcommand{\eenum}{\end{enumerate}}

\renewcommand{\empty}{\emptyset}
\renewcommand{\^}[1]{\hat{#1}}

\newtheorem{definition}{Definition}[section]
\newtheorem{lemma}[definition]{Lemma}
\newtheorem{theorem}[definition]{Theorem}
\newtheorem{conjecture}[definition]{Conjecture}
\newtheorem{proposition}[definition]{Proposition}
\newtheorem{corollary}[definition]{Corollary}

\newtheorem{problem}[definition]{Problem}

\newtheorem{claim}[definition]{Claim}
\newtheorem{subclaim}[definition]{SubClaim}
\newtheorem{fact}[definition]{Fact}

\newtheorem{remark}[definition]{Remark}
\newenvironment{proof}{{\bf Proof}: }{\hfill \qed \vspace{0.25in}}
\begin{document}
\title{On Hanf numbers of the infinitary order property\thanks{The authors   thank
the United States - Israel Binational Science foundation for supporting this research 
as well as the Mathematics
department of Rutgers University where part of this research was carried out.
}}

\author{ Rami Grossberg 
\\Department of Mathematics
\\Carnegie Mellon University
\\Pittsburgh,  PA 15213
\and Saharon Shelah\thanks{This paper replaces item \# 259 from Shelah's list of publications.}
\\Institute of Mathematics
\\The Hebrew University of Jerusalem
\\Jerusalem,  91094  ISRAEL\\
\&  \\Department of Mathematics
\\Rutgers University
\\New Brunswick, NJ 08902}
\date {\today}
\maketitle 
\newpage

\begin{abstract} We study several cardinal, and ordinal--valued functions that are  
relatives of Hanf numbers.
Let $\kappa$  be an infinite cardinal, and let $T\ci L_{\kappa^+,\omega}$ be a theory 
of cardinality  $\leq\kappa$, and let $\gamma$  be an ordinal $\geq\kappa^+$. 
For example we look at  \begin{enumerate}
\item
$\mu_{T}^*(\gamma,\kappa):=\min\{\mu^* : \forall\phi\in L_{\infinity,\omega}$, with
$rk(\phi)<\gamma$, 
\underline{if} $T$ has the $(\phi,\mu^*)$-order property \underline{then}
  there exists a formula 
$\phi^{'}(x;y)\in L_{\kappa^+,\omega}$, such that  for every $\chi\geq\kappa,\;\;$ 
$T$ has the $(\phi^{'},\chi)$-order property \}.
\item $\mu^*(\gamma,\kappa):=\sup \{ \mu_{T}^*(\gamma,\kappa) \; |
\;T\in L_{\kappa^+,\omega} \}$.
\end{enumerate}

We discuss several other related functions,  sample results are:

\begin{itemize}
\item It turns out that \underline {if}  $T$ has the $(\phi , \mu^*(\gamma,\kappa))$-order
propery for some $\phi\in L_{\infinity,\omega}$, with
$rk(\phi)<\gamma$ \underline {then}  for every $\chi > \kappa$ we have that $I(\chi,T)=
2^{\chi }$ holds.
\item For every $\kappa$ and $\;\gamma$ as above there exists an ordinal
$\delta^*(\gamma,\kappa)$  such that  
$\mu^*(\gamma,\kappa)=\beth_{\delta^*(\gamma,\kappa)}$,  
\item $\delta^*(\gamma,\kappa)\leq(|\gamma|^{\kappa})^+$,  
\item for  $\kappa$ with uncountable cofinality, we have that 
$\delta^*(\gamma,\kappa)>|\gamma|^{\kappa}$ and 
\item the ordinal $\delta^*(\gamma,\kappa)$ is bounded below by 
the Galvin--Hajnal rank of
a reduced product.
\end{itemize}

For many cardinalities we have better bounds, some of the bounds obtained using  
Shelah's {\em PCF theory}.  The function  $\mu^*(\gamma,\kappa)$ is used to compute
bounds to 
the values of the function  $\overline {\mu}(\lambda,\kappa)$ we studied in a previous paper.
\end{abstract}
\newpage

\section{Introduction} 

Let $\chi$ be an infinite cardinality, and suppose that $T\ci L_{\chi,\omega}$  
(notice that
when $\chi=\omega$ we are dealing  with first-order theories).  

The fundamental meta-problem in the area of classification theory  can be stated
as:

\begin{problem} \label {basic.prob} What is the structure of $Mod(T)$?
\end{problem}

A more precise (and concrete) test-question is:

\begin{problem} \label {spectrum.prob} What are the possible functions  
$I(\cdot , T): Card \rightarrow Card$?  (where $I(\lambda ,T)$  stands for
the number of isomorphism types of models for $T$ of cardinality $\lambda$)
\end{problem}

A much more precise (and a very difficult) particular case of \ref{spectrum.prob} is the 
following 
\begin{conjecture} {\rm (Shelah about 1976)}  Let  $\psi \in L_{\omega_1 ,\omega}$ be given.
If  there exists a cardinality $\mu > \beth _{\omega_1}$ such that $I(\mu ,\psi )=1$
 then for every $\mu > \beth _{\omega_1}$,  $I(\mu ,\psi )=1$ holds.
\end{conjecture}

A possible approach to Problem \ref{basic.prob} and its relatives, is to try to imitate
Classification theory for elementary classes (see \cite{bible}).  Namely it would be desirable
to find properties parallel to {\em stability}, {\em superstability} etc.  Much work 
has been done in the last 25 years  (see for example -- \cite{Sh 48},\cite{Sh 87a},
\cite{Sh 87b},\cite{Sh 88},\cite{Sh 300},\cite{MaSh 285},\cite{GrSh238}, or \cite{Sh 299}
for a general survey).  In this article we concentrate on
dealing with the parallel (for infinitary languages) to instability.  The following can 
be viewed as a definition of stability for first-order theories:

\begin{fact} \label {unstable}{\rm (\cite{Sh16})}  Let  $T$ be a complete first-order theory.  The following are 
equivalent:
\begin{enumerate}
\item  $T$ is unstable
\item  There exist a formula $\phi (x;y)\in L(T)$,  a model
$M$ for $T$, and a set $\{a_n : n<\omega \}\; \subseteq M$  such that 
$l(x)=l(y)=l(a_n)$ for every $ n<\omega$,  and  for all $n,k<\omega$  we have 
$n<k \Longleftrightarrow  M\models \phi [a_n;a_k].$\\
\end{enumerate}

\end{fact}

Condition 2 in Fact \ref{unstable} is sometimes called the \underline {order-property}.
One of the most important properties of unstable theories is the following:

\begin{fact}\label{many}{\rm \cite{Sh12}}  Let  $T$ be a complete first-order theory.  If $T$ is
unstable  then  for every $\mu >|T|$  we have that $I(\mu , T)=2^{\mu}$.
\end{fact}

An inspection of  the proof of \ref{many} shows that the hypothesis that $T$ is 
a complete first-order unstable theory could be replaced by the following
property:

\parbox{0.0\linewidth} \hfill $(*)_{{}_{T}}\;\;$ \hfill
\parbox{0.88\linewidth}   {There is an expansion $L'$ of  
$L(T)$ with built-in Skolem functions and an $L'$-structure
$M$,  a ``formula'' $\phi(x;y)$ , and there exists $I:=\{a_i : i<\omega \}\subseteq M$ 
a sequence of $L'$-indiscernibles
such that $l(x)=l(y)=l(a_n)$ for every $ n<\omega$, $M$ is the Skolem Hull of $I$, 
$M\restriction L(T)\models T$,
  and for every $n,m < \omega$ we have\\
$n<m\Longleftrightarrow M\models \phi[a_n;a_m].$}
By ``formula'' we mean 
that $\phi$ is  in any logic (over the vocabulary $L'$) such that  $\phi$ is 
preserved by isomorphisms of $L(T)$-structures.
 \\

The condition in Fact \ref{unstable} seems to be a natural candidate for a
definition of instability for infinitary logics.  Since the compactness
theorem fails even for $L_{\omega_1,\omega }$ the next definition
is the replacement of the (the first-order) order-property.

\begin{definition}  Let  $T\ci L_{\chi,\omega}$,  
$\phi(x;y)\in L_{\infty,\omega}$,  and  let  $\mu$  be a cardinality.  \begin{enumerate}
\item We say that
 \underline 
{$M$ has the $(\phi,\mu)$-order property} iff   there exists
$\{a_i : i<\mu\}\ci M$  such that $l(x)=l(y)=l(a_i)<\omega $, and for every $ i,j<\mu$ we have 
 $i<j\iff M\models\phi[a_i;a_j]. $
\item  \underline 
{$T$ has the $(\phi,\mu)$-order property} iff   there exists  $M\models T$  such that 
$M$  has the $(\phi,\mu)$-order property.
\item \underline 
{$T$ has the $(\phi,\infinity)$-order property} iff   for every $\mu$,  $T$ has the 
$(\phi,\mu)$-order property.
\item
Let $\lambda$ and $\mu$ be cardinalities, we say that \\
\underline 
{$T$ has the $(L_{\lambda,\omega},\mu)$-order property} iff there exists 
$\phi\in L_{\lambda,\omega}$  such that $T$ has the 
$(\phi,\mu)$-order property.

\end{enumerate}
\end{definition}

\begin{remark}\label{ordermany} {\rm 
\begin{enumerate}
\item In light of the last definition, 
Fact \ref{unstable} can be restated as (for first-order complete $T$):
  $T$ is unstable iff $T$ has the $(L_{\omega ,\omega },\aleph_0 )$-order property.
\item
It is not difficult to see  (using \cite{Mo}, see \ref{mufact} below) that 
the following implication is true: If $T$ has the 
$(L_{\lambda ^+,\omega},\infinity )$-order property 
then  $(*)_{T}$ holds.
 
\end{enumerate}
}
\end{remark}

The natural question to ask in this context is: Given a theory  $T$ and a cardinality $\mu$,
 does $T$  have the  $(L_{\lambda ^+,\omega},\mu )$-order property?  
The main object of study in \cite{GrSh} was the function ${\overline \mu }(\lambda ,\kappa )$.
The following 
$\mu^*(\lambda,\kappa)$ is a relative of  ${\overline \mu }(\lambda ,\kappa )$ 
from \cite{GrSh}.

\begin{definition}\label{mu*}  Let  $\kappa\leq\lambda$.\begin{enumerate}
\item Let $\psi\in L_{\kappa^+,\omega}$, 
$\mu_{\psi}^*(\lambda,\kappa):=\min\{\mu^* : \forall\phi\in L_{\lambda^+,\omega}$ 
\underline{if} $\psi$ has the $(\phi,\mu^*)$-order property, \underline{then}
 $\exists
\phi^{'}(x;y)\in L_{\kappa^+,\omega}$, such that  $\psi$ 
has the $(\phi^{'},\infinity)$-order property \}.
\item $\mu^*(\lambda,\kappa):=\sup\{\mu_{\psi}^*(\lambda,\kappa)\;|
\; \psi(x;y)\in L_{\kappa^+,\omega}\}$.
\end{enumerate}
\end{definition}

\begin{remark}{\rm
The idea behind Definition \ref{mu*} is that when  $\psi$ has the \\
$(L_{\lambda^+,\omega},\mu^*(\lambda,\kappa))$-order property then (by Remark \ref{ordermany}
and $(*)_{\psi}$)  for every  $\chi >\kappa\;\; I(\chi,\psi )=2^{\chi}$.
}
\end{remark}

Already in \cite{Sh16}  Shelah realized the importance of the above concept (it did not
appear there explicitly. Only in \cite{Gr} (see \cite{GrSh}) we realized the importance of 
functions of this kind).  The previous definition is a generalization of one of our
definitions from \cite{GrSh}, see Definition \ref{mu*}.  
Shelah's fundamental result from \cite{Sh16} 
can be restated as:

\begin{fact}\label{mufact}
{\em (\cite{Sh16})}  For every  $\kappa\leq\lambda$, we have $\mu^*(\lambda,\kappa)\leq
\mu_0(\lambda,1)$\footnote{$\mu_0(\lambda,\lambda)$ is the usual Morley number to 
be introduced in Definition \ref{Morleyno} below. It is easy to see that 
$\mu_0(\lambda,\lambda)=\mu_0(\lambda,1)$}.
\end{fact}


 Let us mention here the following dramatic improvement 
(for  $\kappa=\aleph_0$) of \ref{mufact}:

\begin{theorem}{\rm (\cite{GrSh})}
  For every  $\lambda\geq \aleph_0$, we have $\mu^*(\lambda,\aleph_0)\leq
\beth_{\lambda^+}$.
\end{theorem}

It turns out that even for first-order theories the above question is interesting 
(for
$\kappa  =|L|=\aleph_0$,  $T$ is a complete first-order theory in $L$, we could ask 
what is
an upper bound of $\mu_{T}^*(\lambda,\aleph_0)$?).
Since there are cases 
when $T$ is  stable  (i.e. there is no first-order formula defining an 
$\omega$-sequence
in a model of $T$)  but still $T$ has a hidden instability (like in the case
of stable theories with the omitting-types order property).

Namely there is a natural class of examples of theories that do not have a  first-order 
formula exemplifying the
 order-property but do have an infinitary order property.  
Any stable first order theory that has the omitting types order-property has the
 $(L_{\omega_1 ,\omega },\infinity )$-order property but not the  
$(L_{\omega ,\omega },\aleph_0 )$-order property (see \cite{Sh 200}).

Already from Morley's omitting-types theorem it follows that given  
$T$ and $\phi$ as above there exists  $\mu:=\mu(T,\phi)$  such that 
if  $T$ has the $(\phi,\mu)$-order property then $\forall\lambda\geq\chi$, $T$ has the
$(\phi,\lambda)$-order property.  The bound obtained from repeating the argument 
in the proof of
Morley's omitting types theorem (see \cite{Sh16}) is: 
$\mu(T,\phi)\leq\max\{Hanf(T),Hanf(\phi)\}$.
 Where $Hanf(T)$ and $Hanf(\phi )$ are the Hanf numbers of $T$ and the logic 
containing $\phi$ (respectively).

Let $\chi>\aleph_0$  ($T$ still may be first-order).  
Our object is to find upper bounds on $\mu$.  
It turns out that for
 $\phi\in L_{\infinity,\omega}-L_{\chi,\omega}$ there is a  cardinality  
$\mu^*:=\mu^*(T,\phi)$\footnote{
The surprise is that often  $\mu^*(T,\phi)$ is much smaller than $\mu (T,\phi )$.}, 
such that the following implication holds:  If $T$ has the 
$(\phi,\mu^*)$-order property then there exists a formula $\phi'\in L_{\chi,\omega}$ 
(it is a collapse of $\phi$)  such that $T$ has the $(\phi ',\lambda)$-order property
for every $\lambda\geq\chi$.

In this paper we present a systematic study of several
cardinal and ordinal valued functions related to the infinitary order
property. This  is a continuation of \cite{GrSh},  we deal with similar
problems and improve many of our results.  This is achieved via a 
generalization of the original
problem (dealing with new cases)  while obtaining often
better estimates to our earlier bounds.  
The reader is not expected  to be familiar
with \cite{GrSh}.\\

$\bf {Notation}$:  Everything is standard.   Often  $x$, $y$, and $z$  will 
denote free variables or finite
sequences of variables,  when $x$  is a sequence  $l(x)$  denote its length.  
It should be clear from the context whether we deal with variables or
sequences of variables.  $L$ will denote a similarity type 
(also known as-language or signature),  $\Delta$
will stand for  a  set of  $L$ formulas.  $M$  and $N$  will stand for  
$L$ - structures, $|M|$ the universe of the structure $M$, $\|M\|$  
the cardinality of the universe of $M$.
Given a fixed structure  $M$,  subsets of its universe will be denoted by 
$A$, $B$, $C$, and $D$.
So when we write  $A \ci M$ we really mean that $A \ci |M|$,  while $N \ci M $ 
 stands for  
``$N$  is a submodel of  $M$''.  Let $M$  be  a structure.  By  $a \in M$ we mean  
$a \in |M|$,   
when $a$  is a finite sequence of elements then  $a \in M$ stands for ``all the elements
of the sequence  $a$  are elements of $|M|$''.  For cardinalities $\kappa\leq\lambda$,
 let $S_{<\kappa}(\lambda):=\{X\ci \lambda \; :\; |X|<\kappa \}$. 
When  $T$ is a first-order theory, 
$\Gamma $  denotes a set of  $T$-types over the empty set (not necessarily complete types).
$EC(T,\Gamma):=\{M\;:\; M\models T,\;\forall p\in\Gamma\mbox{ $M$ omits the type $p$} \}$.
When $T$ is first-order, $L\subseteq L(T)$, and $\Gamma$ is a set of $T$-types by  
$PC(T,\Gamma , L)$ we denote the following 
$\{ M\restriction L\;:\;M\models T,\;\forall p\in\Gamma\mbox{ $M$ omits the type $p$} \}$; 
namely  $EC(T,\Gamma)=PC(T,\Gamma ,L(T))$.
$\lambda,\mu,\kappa$, and $\chi$ will stand for infinite cardinalities;  
$\alpha,\beta,\gamma,\delta,
\zeta$, and $\xi$ are ordinals.
References of the form ``Theorem IV 3.12'' are to \cite{bible}.
For $\phi\in L_{\infinity ,\omega}$,  let $Sub(\phi)$ be the set of subformulas of 
$\phi$, now let
\[ rk(\phi ):=\left\{\begin{array}{ll}
$0$ & \mbox{if $\phi$ is atomic} \\
Sup\{rk(\chi )+1 \;:\; \chi \in Sub(\phi)\} & \mbox{otherwise.} \end{array} 
\right. \]
\newpage
\section{Review} 

C.C. Chang in  \cite{Ch}  made the following fundamental observation:

\begin{fact}\label{Ch}  Let  $\kappa$ be an infinite cardinality,  and let 
$L$ be a similarity type of cardinality no more than $\kappa$.  
Given  $\psi\in L_{\kappa^+ , \omega}$, there exist  a similarity type 
$L'\supseteq L$, a 
first-order theory  $T$ in $L'$, and a set of $T$-types $\Gamma$  
(all three of cardinality less 
than
or equal to $\kappa$) such that
$Mod(\psi)=PC(T,\Gamma , L)$.
\end{fact}

Namely instead of studying $Mod(\psi)$ directly for an infinitary theory $\psi$ it is enough
to look at a  class of reducts of models of a first-order theory that omits a set
of types.

W. Hanf and M. Morley \cite{Mo}, 
recognized the importance of the following concept:

\begin{definition}\label{Morleyno}
Let  $T$  be a first-order theory,  and let  $\Gamma$ be a set of $T$-types. \\ 
\underline {The Morley number}\footnote{Some authors call this the Hanf number of $T$ and the 
$\Gamma$}
of $T$ and $\Gamma$, is the following:\begin{enumerate}
\item $\mu_0(T,\Gamma):=\min \{ \mu \;:\;\exists M\in EC(T,\Gamma) \;\|M\|\geq\mu\Rightarrow
\forall\chi\geq|T|\;\;\exists N\in EC(T,\Gamma) $ of cardinality $\geq\chi$\}.
\item Let  $\lambda,\;\kappa$ be cardinalities.\\
$\mu_0(\lambda,\kappa):=\sup\{\mu_0(T,\Gamma)\; :\; |T|\leq\lambda,  \; \Gamma$ a set of $T$-types
 of 
cardinality  $\leq\kappa$\}\footnote{We hope that the reader is not bothered by this abuse 
of notation.  We are using the same letter $\mu_0$  to denote entirely different 
(but related) functions.
They can be distinguished by the type of the arguments they take.}.
\end{enumerate}
\end{definition}

Morley (among other things)  have shown   that $\mu_0(\aleph_0,\aleph_0)=\beth_{\omega_1}$.
His most general result is stated as Theorem \ref{Morleydelta} below.  Shelah in \cite{Sh78}
have dealt with what can be viewed to be an interpolant of $\mu_0(T,\Gamma)$ and
 $\mu_0(\lambda,\kappa)$:
\[\mu_0(T,\kappa):=\sup\{\mu_0(T,\Gamma)\;:\; \Gamma \mbox{ a set of $T$-types, }\; |\Gamma|\leq\kappa \}.\]

It is not difficult to conclude from the proof of Morley's categoricity theorem that 
when $T$ is a countable  and 
$\aleph_0$-stable theory then $\mu_0(T,\cdot)\leq\aleph_1$.  Shelah in \cite{Sh78}
 studied the effect  that stability of $T$ has  on the upper bounds on $\mu_0(T,\kappa)$.
This work was continued about ten years later by Hrushovski and Shelah in \cite{HrSh}.

In this paper, since our main goal is the study of unstable theories 
(or theories that are not stable in a weak sense) we will ignore the effect that  
the stability of $T$ may have  on
the function $\mu_0(T,\kappa)$.

The modern era in the study of Hanf numbers began with the paper of 
Barwise and Kunen \cite{BaKu}.  They studied systematically the relationship between
the function $\mu_0$ and the first ordinal that exemplify the undefinability of well
ordering in classes of models that omit a set of types. 
 Below we introduce an ordinal-valued
function $\delta_0(\lambda , \kappa )$ that turns out to be related to 
$\mu_0 (\lambda ,\kappa )$
in a nice way.

\begin{definition}
Let  $\lambda$  and  $\kappa$  be infinite cardinalities,  $T$ varies over consistent
first-order theories such that $L(T)\supseteq\{P,<\}$ when $P$ is a unary 
predicate and $T\vdash$``$<$ linearly orders $P$''.\\
 $\delta_0(\lambda,\kappa):=\min\{\delta : |T|\leq\lambda, 
\Gamma \mbox{ a set of $T$-types, } 
|\Gamma|\leq\kappa$  \underline{if} for every $\delta^{'} <\delta$ \\
there exists 
$M\in EC(T,\Gamma)$
such that $otp(P^M , <^M )\geq\delta^{'} $  \underline{then} there \\exists  
$N\in EC(T,\Gamma)$  s.t. 
$(P^N,<^N)$ is not well ordered \}.

\end{definition}

 The following is a restatement of Morley's ``other'' important theorem:

\begin{fact}\label{Morleydelta}
{\em (Theorem VII 5.5)}  $\mu_0(\lambda,\kappa)=\beth _{\delta_0(\lambda,\kappa)}$.
\end{fact}

The following ordinal and cardinal-valued  functions are from \S 4 of \cite{GrSh}:

\begin{definition}\label{delta1}
Suppose $T$ is a first-order theory such that $L(T)$  contains
 $\{<,P\}$ and
\[ T\vdash [<\mbox{ is a linear order}]\;\wedge \;[\;<\restriction P 
\mbox{ is a linear order on the unary predicate } P ]. \]  \begin{enumerate}
\item $\delta_1(\theta,\lambda,\kappa):=\min\{\delta : \;\Gamma \mbox{ a set of $T$-types, }
\;|\Gamma|\leq \lambda, |T|\leq\kappa\;$ \\\
\underline{if} $\forall\delta^{'}<\delta\;\;\exists M\in EC(T,\Gamma)\;\;$ with 
$otp(P^M,<^M)\in On\cap\theta^+$  and  $otp(M-P,<^M)\geq\delta^{'}$, \underline{then} 
$\exists N\in EC(T,\Gamma) $  s.t.  $otp(P^N,<^N)\in On\cap\kappa^+$  and 
$(N-P^N,<^N)$ is not well ordered \}.
\item $\mu_1(\theta,\lambda,\kappa):=\min\{\mu : 
|\Gamma|\leq \lambda, |T|\leq\kappa\;
\;\Gamma \mbox{ a set of $T$-types, }\;
$ \underline{if} $\exists M\in EC(T,\Gamma)\;\;
 \|M\|\geq\mu $  with $otp(P^M,<^M)\in On\cap\theta^+$ \underline{then}
 for every $\chi\geq\kappa\;\;\exists
N\in EC(T,\Gamma) $  of cardinality at least $\chi$ such that
$otp(P^N,<^N)\in On\cap\kappa^+$ \}.
\item  When $\theta=\lambda$ we will omit the first parameter-$\theta$
\end{enumerate}
\end{definition}

We will prove the following equality:

\begin{theorem}\label{ineq}
For every $\kappa\leq\lambda\leq\theta$  we have  
$\mu_1(\theta,\lambda,\kappa)=\beth_{\delta_1(\theta,\lambda,\kappa)}$.
\end{theorem}
 
From now on we concentrate on the case that $\theta=\lambda$ and we will work
with the functions $\delta_1(\lambda,\kappa)$ and $\mu_1(\lambda,\kappa)$.  The arguments for 
the functions with three parameters are essentially similar (they require an additional technical
effort, but require no new ideas).  Note that by \cite{Ch} working with two
parameter functions is sufficient for $L_{\lambda^+,\omega}$.
The new point is that we are able to show that  
$\mu_1(\lambda,\kappa)\geq\beth_{\delta_1(\lambda,\kappa)}$.  
The proof of Theorem \ref{ineq} 
is similar to that of Theorem VII 5.5, we skip its proof, 
since later we will prove a related theorem (Th. \ref{eq2}) whose proof 
is similar (and little harder).

The next proposition provides us with a lower bound for
$\delta_1(\lambda,\kappa)$, it follows immediately from the definitions (in Theorem
 \ref{lowerbound} we show a better lower bound).

\begin{proposition}
For $\lambda\geq\kappa$ we have that 
$\delta_1(\lambda,\kappa)\geq\delta_1(\kappa,\kappa)\geq
\delta_0(\kappa,\kappa)=\delta_0(\kappa,1)$.
\end{proposition}

In the following proposition the connection between the last definition and 
the order property is clarified.

\begin{theorem}\label{mu1=mu*}  Let  $\kappa\leq \lambda$, be cardinalities. 
$\beth_{\lambda^+}\leq\mu^*(\lambda,\kappa)\leq\mu_1(\lambda,\kappa)$.  

\end{theorem}
{\em Proof.}  First we show that $\mu^*(\lambda,\kappa)\leq\mu_1(\lambda,\kappa)$.
Let $\psi\in L_{\kappa^+,\omega}$, and  
$\phi(x;y)\in L_{\lambda^+,\omega}$  
be given.  Suppose  $\psi$  has the $(\phi,\mu_1(\lambda ,\kappa))$-order property
we need to find a formula  $\phi'\in L_{\kappa^+,\omega}$ 
such that $\psi$ has the
$(\phi',\infinity)$-order property.

By Fact \ref{Ch} there exists a first-order theory $T$ in a similarity type $L(T)$
that extends $L$, and there is a set $\Gamma$ of $T$-types of cardinality $\leq\kappa$
such that $PC(T,\Gamma,L)=Mod(\psi)$.
 By following the inductive definition of the formula $\phi$
we may identify $\phi$ with a function $f$ from the set $P$
into the set $ L\union \{\wedge,\neg,(,),=\}\union \{x_i\;:\;i<\kappa\}$.

Let $\chi$ be a regular large enough such that 
\[ \{ L_{\lambda^+,\omega},L_{\kappa^+,\omega},T,\Gamma,L,\phi,f,P,\psi,\lambda^+,
\mu_1(\lambda ,\kappa),\delta_1(\lambda,\kappa) \}\union\mu_1(\lambda ,\kappa)
\subseteq H(\chi). \]
In addition we require that the structure $\langle H(\chi),\in\rangle$
 reflects all the relevant
properties of the above sets.  Let $P$ be the rank of the formula $\phi$,
 note that it is an
ordinal less than $\lambda^+$.   
Let ${\mathfrak A'} \prec \langle H(\chi),\in,\ldots\rangle$
of cardinality $\mu_1(\lambda ,\kappa)$ such that $\mu_1(\lambda ,\kappa)^{\mathfrak A'}=
\mu_1(\lambda ,\kappa)$ (so $\mu_1(\lambda,\kappa)+1\subseteq \mathfrak A'$),
 fix a bijection $G$ from $\mu_1(\lambda ,\kappa)$  onto the
universe of  ${\mathfrak A'}$, and let  ${\mathfrak A}:=\langle  {\mathfrak A'},G\rangle$.
By the definition of
$\mu_1(\lambda,\kappa)$, for every $\chi\geq\kappa$ there exists ${\mathfrak B}_{\chi}\equiv 
{\mathfrak A}$ 
of cardinality $\chi$ such that
 ${\mathfrak B_{\chi}}$ omits the  types from $\Gamma$, 
$\kappa^{ {\mathfrak B_{\chi}}}=\kappa$, and  $P$
 is an ordinal less than $\kappa^+$ (just apply
the Mostowski collapse on $\mathfrak B_{\chi}$).
  Using $P^{ {\mathfrak B_{\chi}}}$,
 and $f^{ {\mathfrak B_{\chi}}}$ we
know (in ${\mathfrak B_{\chi}}$) that $\phi^{ {\mathfrak B_{\chi}}}
\in L_{\kappa^+,\omega}^{ {\mathfrak B_{\chi}}}$,
 but since
$\kappa^{ {\mathfrak B_{\chi}}}=\kappa$ we have that $\phi^{
 {\mathfrak B_{\chi}}}
\in L_{\kappa^+,\omega}$ is
a formula as required in the definition of $\mu^*(\lambda,\kappa)$.

To see that  $\mu^*(\lambda,\kappa)\geq\beth_{\lambda^+}$: It is enough to show that for
every $\alpha<{\lambda^+}$ there exist a sentence $\psi_{\alpha}\in L_{\kappa^+,
\omega}$, and a formula $\phi_{\alpha}\in L_{\lambda^+,\omega}$ such that 
$\psi_{\alpha}$ has
the $(\phi_{\alpha},\beth_{\alpha})$-order property and $\psi_{\alpha}$ 
does not have the
 ($L_{\kappa^+,\omega},\infinity)$-order property.

Before proving this we will introduce several definitions.\\

$\bf {Notation}$:
The sentence $\psi_{\alpha}$ will be defined as a the theory of a well founded tree.
We  deal with well-founded trees whose vertices are
 decreasing sequences of ordinals, the root of the 
tree ${\cal T}$ is denoted by $rt({\cal T})$,
for an element $x\in {\cal T}$ let $Suc_{\cal T}(x)$ 
stand for the set of immediate succssesors of $x$,
 and
${\cal T}[x]$ stands for the subtree of ${\cal T}$
 consisting of the elements that are greater  than or equal 
to $x$.

\begin{definition} 
Let ${\cal T}$ be a well founded tree.  
\begin{enumerate}
\item For $x\in {\cal T}$ let $Dp_{\cal {\cal T}}(x)=\beta$ the
 \underline {depth of $x$ in ${\cal T}$} defined
by induction on $\beta$:
\begin{enumerate}
\item if $Suc_{\cal T}(x)=\empty$ then $Dp_{\cal T}(x)=0$.
\item if for every  $y\in Suc_{\cal T}(x)$  we have $Dp_{\cal T}(y)<\beta$, and for every 
$\gamma<\beta$ there exists  $z\in Suc_{\cal T}(x)$ of  such that $Dp_{\cal T}(z)\geq\gamma$ 
then $Dp_{\cal T}(x)=\beta$.
\end{enumerate}
\item The \underline {depth of ${\cal T}$} is $Dp({\cal T}):=
\sup\{\; Dp_{\cal T}(x)\;:\;x\in {\cal T}\}$.
\end{enumerate}
\end{definition}

\begin{proposition}\label{DT=DT}
Let ${\cal T}$ be a well-founded tree,  $Dp({\cal T})=Dp_{\cal T}(rt({\cal T}))$.
\end{proposition}
{\em Proof.}Trivial.              \hfill $\square_{\ref{DT=DT}}$

\begin{claim}\label{claim}
For every  $\alpha$ there exists a well-founded tree ${\cal T}_{\alpha}$ of depth 
$\alpha$ such that
$\|{\cal T}_{\alpha}\|\leq |\alpha|+\aleph_0$.
\end{claim}
{\em Proof.}
By induction on $\alpha$:

For $\alpha=0$; Simply let ${\cal T}_0:=\langle\rangle$. 

For $\alpha=\beta+1$; Suppose ${\cal T}_{\beta}$ is a tree of depth $\beta$. \\
Let ${\cal T}_{\alpha}:=\{\langle\rangle\}\union\{\langle\rangle \^ \;\eta \;:\;\eta
\in {\cal T}_{\beta } \}$.  
The order on ${\cal T}_{\alpha}$ is the obvious.

For $\alpha$ a limit ordinal; By the induction hypothesis let
 $\{{\cal T}_{\beta}\;:\;\beta<\alpha \}$ be pairwise disjoint trees, each
of depth $\beta$. 
\\ Define ${\cal T}_{\alpha}$ to be the tree 
$\{\langle\rangle\} \union\{\langle\rangle \^ \; \eta \;:\;\eta\in
 {\cal T}_{\beta}, \beta<\alpha\}.$
\hfill $\square_{\ref{claim}}$

\begin{definition}
\begin{enumerate}
\item Let ${\cal T}_1$, ${\cal T}_2$ be well-founded trees, and let $\alpha$ be an ordinal.  
By induction on
$\alpha$ define when ${\cal T}_1\approx_{\alpha}{\cal T}_2$:
\begin{enumerate}
\item For $\alpha=0$, always ${\cal T}_1\approx_{\alpha}{\cal T}_2$.
\item For $\alpha\neq 0$, if for every $\beta<\alpha$ and for every 
$x_1\in Suc_{{\cal T}_1}(rt({\cal T}_1))$ 
there exists $x_2\in Suc_{{\cal T}_2}(rt({\cal T}_2))$ such that 
${\cal T}_1[x_1]\approx_{\beta}{\cal T}_2[x_2]$,
 and for every
 $x_2\in Suc_{{\cal T}_2}(rt({\cal T}_2))$ 
there exists $x_1\in Suc_{{\cal T}_1}(rt({\cal T}_1))$ such that 
${\cal T}_2[x_2]\approx_{\beta}{\cal T}_1[x_1]$.
 \end{enumerate}
\item A tree ${\cal T}$ is called \underline {simple} iff there are no distinct 
$x_1,x_2\in Suc_{\cal T}(rt({\cal T}))$ such that $Dp_{\cal T}(x_1)=Dp_{\cal T}(x_2)$ and 
${\cal T}[x_1]\approx_{Dp_{\cal T}(x_1)}{\cal T}[x_2]$.
\end{enumerate}
\end{definition}

\begin{proposition}\label{prop1}
Let  ${\cal T}_1,{\cal T}_2$ be trees, and let  $\alpha$ be an ordinal.  
If ${\cal T}_1\approx_{\alpha}{\cal T}_2$ then 
one of the following conditions holds:
\begin{enumerate}
\item $Dp({\cal T}_1)=Dp({\cal T}_2)$, or
\item $Dp({\cal T}_1)\geq\alpha$ and $Dp({\cal T}_2)\geq\alpha$.
\end{enumerate}
\end{proposition}
{\em Proof.} Easy, by induction on $\alpha$.
\hfill $\square_{\ref{prop1}}$

\begin{claim}\label{existenceoftrees}
For every ordinal $\alpha$ there exists a family of simple trees $\{{\cal T}_i \; : \; 
i<\beth_{\alpha} \}$, such that for every $i<\beth_{\alpha}$
\begin{enumerate}
\item $\| {\cal T}_i \|\leq\beth_{\alpha}$,
\item $i\neq j\Rightarrow {\cal T}_i\not\approx_{\omega+\alpha}{\cal T}_j$,
\item $Dp({\cal T}_i)=\omega+\alpha+1$.
\end{enumerate}
\end{claim}

{\em Proof.}  By induction on $\alpha$:

For $\alpha=1$; First construct $\aleph_0$ simple trees $\{{\cal T}_n\subseteq^
{\;\omega >}\omega \; : \; n<\omega \}$ 
such that ${\cal T}_0:=\{\langle\rangle\}, \; {\cal T}_{n+1}:=\{\langle n+1\rangle\}\union {\cal T}_n$ 
when the order is an extention
of the order on ${\cal T}_n$; $\langle\rangle$ is the root, and $\langle n+1\rangle$ is a new 
immediate successor of the root incomparable with the elements of 
$Suc_{{\cal T}_n}(\langle\rangle)$.
Now for every $A\subset\omega$ let 
\[{\cal T}_A:=\{\langle \omega\rangle \^ \; \eta \; : \; \eta\in {\cal T}_n,\;n\in A\}
\union \{\nu\in {\cal T}_k \; : \; k\not\in A \}. \]

The order of ${\cal T}_A$ is defined as follows:  $\langle \omega\rangle$ 
is a new immediate sucessor
of the root, the elements $\langle \omega\rangle \^ \;\eta$ and $\nu$
 are pairwise incomparable
when $\eta\in {\cal T}_n,\;(n\in A)$ and $\nu\in {\cal T}_k,\; (k\not\in A)$, and we require that \\
$\langle \omega\rangle \^ \; \eta_1<\langle \omega\rangle \^ \;\eta_2\Leftrightarrow
\eta_1<_{{\cal T}_n}\eta_2$. \\
 In order to see that $A\neq B\subseteq\omega\Rightarrow {\cal T}_A\not
\approx_{\omega+\alpha} {\cal T}_B$: W.l.o.g. we may assume that $\exists n\in A-B$. 
Since ${\cal T}_A[\langle \omega,n\rangle]\not\approx_{\omega}{\cal T}_B[\nu]$ for any 
$\nu\in Suc_{{\cal T}_B}(\langle \omega\rangle)$  (this is because ${\cal T}_k$ and ${\cal T}_n$ are 
inequivalent for $k\neq n$).

For $\alpha\neq 1$; By the inductive hypothesis let $\{\{{\cal T}^{\beta}_i \; : \; i<
\beth_{\beta}\}
 \; :\; \beta < \alpha \}$ be disjoint trees satisfying the statement of the Theorem.
Denote by $S$ the  set 
$\{\langle\beta,i\rangle \;:\;i<\beth_{\beta},\; \beta<\alpha \}$. 
 Fix an injective mapping
from $S$ into $On-Sup(\bigcup_{\beta,i}{\cal T}^{\beta}_i) $, 
denote by $\gamma_{\beta,i}$ the image of the pair $\langle\beta,i\rangle$.
For every $\gamma<\alpha$, and for every 
$A\subseteq S$
cardinality $\beth_{\gamma}$ define 
\[{\cal T}_A:=\{\langle\rangle\}\union\{\langle\gamma_{\beta,i}
\rangle \^ \;\eta \; : \;\eta\in {\cal T}^{\beta}_i, \; \langle\beta,i\rangle\in A \}. \]

The order on ${\cal T}_A$ is defined in the natural way:  $\langle\rangle$ is the root, and 
$\langle\beta,i\rangle \^\;\eta_1<_{{\cal T}_A}\langle\beta,i\rangle \^\;\eta_2$ iff 
$\eta_1<_{{\cal T}^{\beta}_i}\eta_2$.  The verification that $Dp({\cal T}_A)=\omega+\alpha+1$ 
 is left to the reader.  Suppose  $A\neq B\subseteq \beth_{\alpha}$ both of cardinality
$\beth_{\beta}$ for some $\beta<\alpha$.  We need to show that 
$ {\cal T}_A\not\approx_{\omega+\alpha}{\cal T}_B$.  W.l.o.g. there exists $\gamma_{\beta,i}\in A-B$.
Since $\langle\gamma_{\beta,i}
\rangle \^ \;\eta \; : \;\eta\in {\cal T}^{\beta}_i\}$ is a subtree of ${\cal T}_A$ and for all
$j<\beth_{\beta}$ we have that 
$j\neq i\Rightarrow {\cal T}^{\beta}_i\not\approx_{\omega+\beta}{\cal T}^{\beta}_j$, from the 
definition of the relation $\approx$, and the fact that  it follows that there is
no ordinal $\epsilon$ such that the tree $\{\epsilon \^\;\eta \;:\;\eta\in {\cal T}^{\beta}_i\}$
does not appear as a subtree of ${\cal T}_B$ it is clear that
 $ {\cal T}_A\not\approx_{\omega+\alpha}{\cal T}_B$.
\hfill $\square_{\ref{existenceoftrees}}$\\

Back to the proof of Theorem \ref{mu1=mu*}:  Let $\alpha<\lambda^+$ be given.  
By Claim \ref{existenceoftrees} there exists a family of 
nonequivalent simple trees  $\{{\cal T}_i \; : \; 
i<\beth_{\alpha} \}$.  By renaming, we may assume that the above trees do not
 contain sequences
of ordinals which are less than $\beth_{\alpha}$.  We define a new tree:
 $M_{\alpha}$ its
set of elements consists of 
\[ \{\langle\rangle\}\union \{\langle i\rangle \;:\; i<\beth_{\alpha} \}\union \{
\langle i,i'\rangle\^{ }\eta \;:\;\eta\in {\cal T}_{i'},\;i'< i,\; i<\beth_{\alpha}\}. \]

We can view  $M_{\alpha}$ as a partially ordered set (by ``being initial segment'').
  We view $M_{\alpha}$ as a model in 
a languge consisting of single function symbol:  
A unary  function $f$ whose interpretation is
the predessor of its argument (if the argument is the root than 
the value is defined to be the root).  Notice that the following
formula (of $L_{\omega_1,\omega}$),  $\bigvee_{k<\omega}[x=f^k(y)]$\footnote{
When $f^k(y)$ stands for $f(\cdots f(y)\cdots)$ $k$-many times.} defines the relation of 
``being an initial segment'' on well founded trees. 

Let $\psi_{\alpha}:=\bigwedge Th_{\omega,\omega}(M_{\alpha})\wedge
(\forall x)\bigvee_{n<\omega}[f^n(x)=\langle\rangle]$.  Namely  $\psi_{\alpha}$ is the 
first-order theory of $M_{\alpha}$ together with the statement that say that every element
is of finite distance from the root.  

Let $\phi_{\alpha}(x,y)$ be the following statement:  $x,y\in Suc(\langle\rangle)$, 
and for every  $x'\in Suc(x)$ there exists $y'\in Suc(y)$ such that ${\cal T}[x']\approx
_{\omega+\alpha +1}{\cal T}[y']$, and there exists  $y'\in Suc(y)$ such that 
for every  $x'\in Suc(x)$ we have that  ${\cal T}[x']\not\approx
_{\omega+\alpha +1}{\cal T}[y']$ holds.

In order to complete the proof of Theorem \ref{mu1=mu*}, it suffices to prove the 
following:

\begin{subclaim}\label{subclaim}
\begin{enumerate}
\item $\phi_{\alpha}(x,y)\in L_{\lambda^+,\omega}$, 
\item $M_{\alpha}$ has the $(\phi_{\alpha},\beth_{\alpha})$-order property,
\item There do not exist a formula $\phi'(x,y)\in L_{\kappa^+,\omega}$ such that 
$\psi_{\alpha}$ has the $(\phi'(x,y),\infinity)$-order property.
\end{enumerate}
\end{subclaim}

{\em Proof.}
\begin{enumerate}

\item  Let  ${\cal T}$ be a well founded tree, and let $\alpha<\lambda^+$ be given,
  it is enough to show by induction on $\alpha$ that
there exists a formula $\chi(x,y)\in L_{\lambda^+,\omega}$ such that  for every
$a,b\in {\cal T}$ we have that  ${\cal T}\models \chi[a,b]$ iff ${\cal T}[a]\approx_{\alpha}{\cal T}[b]$.  
It is easy to check
that the relation $\approx_{\alpha}$ is definable in $ L_{\lambda^+,\omega}$.

\item  Check that for every $i_1,i_2<\beth_{\alpha}$ we have that \\
$i_1<i_2$ iff $M_{\alpha}\models\phi_{\alpha}[\langle i_1\rangle ,\langle i_2\rangle ] $.

\item For the sake of contradiction suppose that there exists a formula 
$\phi'(x,y)\in L_{\kappa^+,\omega}$ such that
 $\psi_{\alpha}$ has the $(\phi',\infinity)$-order property.  
Suppose that $\gamma$ is a limit ordinal $<\kappa^+$ such that the formula
 $\phi'$ has quantifier depth $<\gamma$.  Denote by $\mu$ the cardinality 
$(\beth_{\gamma+1}(|L|))^+$.
Let $N\models\psi_{\alpha}$ be a model of cardinality $\mu$ such that there exists 
$\{a_i\;:\; i<\mu\}$ such that  $l(x)=l(y)=l(a_i)=n<\omega$ and for every $i_1,i_2<\mu$ we have 
$i_1<i_2\Longleftrightarrow N\models\phi'[a_{i_1},a_{i_2}]$ holds.
For every $i<\mu$ fix  $\langle b^i_l\;:\; l<n\rangle=a_i$.  
By the $L_{\omega_1,\omega}$-part of the definition of 
$\psi_{\alpha}$ we have that  $N\models (\forall x)\bigvee_{m<\omega}f^m(x)=f^{m+1}(x)$.
For every $c\in N$ let $m(c):=\min \{m\;:\; N\models f^m(c)=f^{m+1}(c) \}$.  

Since $\mu$ is regular, after renaming we may assume that
 for every $l<n$ there are $k_l<\omega$ such that for every $i<\mu$ we have $m(b^i_l)=k_l$.
By increasing $n$ we may assume that for every $i<\mu$ we have
that $f(b^i_l)\in \{b^i_k \;:\; k<n\}$, and for every $i_1,i_2<\mu$ and every $l_1,l_2<n$
we have \[N\models f(b^{i_1}_{l_1})=b^{i_1}_{l_2}\Longleftrightarrow
N\models f(b^{i_2}_{l_1})=b^{i_2}_{l_2}\;\bigwedge\;
N\models b^{i_1}_{l_1}=b^{i_1}_{l_2}\Longleftrightarrow
N\models b^{i_2}_{l_1}=b^{i_2}_{l_2}.\]
We may also assume that $\langle b_l \;:\; l<n\rangle$ has no repeatition.  
By the $\Delta$-system lemma there exists $s\subseteq n$ such that for every $i_1,i_2<\mu$ and
every $l_1,l_2<n$ we have that  $ b^{i_1}_{l_1}=b^{i_2}_{l_2}\iff l_1=l_2\in s$.

Let $\Phi_{\gamma}$ be the set of $L_{\infinity,\omega}$ 
formulas of quantifier depth $<\gamma$ with
 finitely many free variables.  Clearly $|\Phi_{\gamma}|\leq\beth_{\gamma}(|L|)$ and 
$|{\cal P}(\Phi_{\gamma})|\leq\beth_{\gamma+1}(|L|)<\mu=cf(\mu)$.

Let  $tp_{\gamma}(b_0,\ldots ,b_{m-1};M):=\{\phi({\bar x})\in \Phi_{\gamma}\;:\;
 M\models\phi[b_0,\ldots,b_{m-1}]\}$.  Without loss of generality we may assume that for 
every $i,j<\mu$ we have 
 $tp_{\gamma}(b^i_0,\ldots ,b^i_{m-1};N)=tp_{\gamma}(b^j_0,\ldots ,b^j_{m-1};N)$.
\end{enumerate}
We  obtain a contradiction to the assumption that 
 $\psi_{\alpha}$ has the $(\phi',\infinity)$-order property by proving the following:

\begin{claim}For every $i,j<\mu$ we have
\[ N\models\phi'[b^i_0,\ldots,b^i_{n-1},b^j_0,\ldots,b^j_{n-1}]\iff 
N\models\phi'[b^j_0,\ldots,b^j_{n-1},b^i_0,\ldots,b^i_{n-1}] \]
\end{claim}
{\em Proof.} 
Left to the reader.
\hfill $\square_{\ref{subclaim}}$



\begin{remark}{\rm
In Definition \ref{delta1} we have introduced a third parameter, 
but since it does not add 
anything of substance
(just complicates the notation that may be already little  heavy) we decided to limit
our treatment to the above particular case.  At the end of this section we discuss several 
 generalizations.}
\end{remark}

 Theorem \ref{mu1=mu*} provides 
a  better upper bound than the one in Fact \ref{mufact}:

\begin{corollary}\label{bethfact} For every  $\kappa\leq\lambda$, we have
$\mu^*(\lambda,\kappa)\leq\beth_{\delta_1(\lambda,\kappa)}$.

\end{corollary}

\begin{remark} {\rm Using Facts \ref{mufact} and \ref{bethfact}, one can show that 
for $\kappa\leq\lambda$ we have $\delta_1(\lambda,\kappa)\leq\delta_0(\lambda,\kappa)$.
In \cite{GrSh} we have shown that in many instances the ordinal 
 $\delta_1(\lambda,\kappa)$ 
is much smaller than
$\delta_0(\lambda,\lambda)$  [e.g.  when $\kappa=\aleph_0$,  
 we have
that $\delta_1(\lambda,\kappa)=\lambda^+$,  while for $\lambda=\beth_{\omega_1}$, we have 
$\delta_0(\lambda,\lambda)>2^{\lambda}$. ]}
\end{remark}

In Theorem \ref{ineq} we reduced the  problem of finding estimates for $\mu_1(\cdot,\cdot)$  to 
finding bounds for $\delta_1(\cdot,\cdot)$.  In Fact \ref{covfact}, below
we state an important result from \cite{GrSh}, first we need the following:

\begin{definition} For uncountable $\kappa$,  and $\lambda\geq\kappa$,  denote by
\[\kappa^*:=\left\{ \begin{array}{ll}
\kappa & \mbox{if $cf\kappa=\aleph_0$}\\
\kappa^+ & \mbox{if $cf\kappa>\aleph_0$.}\end{array} \right. \]

$cov(\lambda,\kappa):=\min \{ |F| \;:\; F\ci S_{<\kappa^*}(\lambda), \; 
\forall X\in S_{<\kappa^*}(\lambda) \;\exists\{w_l : l<\omega\}\ci F$,  
such that $ X\ci \bigcup_{l<\omega } w_l \}$.
\end{definition}

Clearly $cov(\lambda,\kappa)\leq\lambda^\kappa$.  
But often $cov(\lambda,\kappa)<\lambda^\kappa$.  In  \cite{Sh g} 
Shelah has a more general function.  Our $cov(\lambda,\kappa)$ is the same
as $cov(\lambda,\kappa^*,\kappa^*,\aleph_1)$ from Definition II 5.2 of \cite{Sh g}.

\begin{fact} \rm {(Theorem 4.4 of \cite{GrSh})}\label{covfact}
Let  $\kappa \leq\lambda$  be infinite cardinalities.\begin{enumerate}
\item if $\kappa=\aleph_0$  then  $\delta_1(\lambda,\kappa)\leq\lambda^+$.
\item if $cf\kappa>\aleph_0$  then 
$\delta_1(\lambda,\kappa)\leq(cov(\lambda,\kappa)+2^{\kappa})^+$.
\item if $cf\kappa=\aleph_0$ then 
$\delta_1(\lambda,\kappa)\leq(cov(\lambda,\kappa)+2^{<\kappa}+\aleph_0)^+$.
\end{enumerate}
\end{fact}

Note that the above innocent looking results are quite powerfull!  
E.g. By a result of \cite{Sh g} (from Chapter XI), \\if 
$(\forall\mu<\chi) [\mu^{\kappa}<\lambda]\bigwedge cf(\chi)=\aleph_0\bigwedge
\chi\leq\lambda < \chi^{\delta+\omega_1}$ then
 we have that 
$cov(\lambda,\kappa)=\lambda$, 
thus $\mu_1(\lambda,\kappa)\leq\beth_{\lambda^+}$, while using Morley's methods we get
only  $\mu_1(\lambda,\kappa)\leq\beth_{(2^{\lambda})^+}$.\\

The following is a generalization of the cardinal--valued function 
we have introduced in 
Definition \ref{mu*}.  Here instead of assuming that $\phi (x;y)$ is an   
$ L_{\lambda^+,\omega}$ formula we look
at all $\phi\in L_{\infinity,\omega}$  with quantifier  depth $\;<\gamma$,
we take into consideration only the 
rank of the formula  $\phi$.

\begin{definition}  Let  $\kappa$  be  an infinite cardinality,  and let $\gamma$  
be an \underline {ordinal}
greater or equal to $\kappa^+$,  $T\in L_{\kappa^+,\omega}$\begin{enumerate}
\item
$\mu_{T}^*(\gamma,\kappa):=\min\{\mu^* : \forall\phi\in L_{\infinity,\omega}$, with
$rk(\phi)<\gamma$, 
\underline{if} $T$ has the $(\phi,\mu^*)$-order property, \underline{then}
 $\exists
\phi'(x;y)\in L_{\kappa^+,\omega}$, such that  $T$
 has the $(\phi',\infinity)$-order property\}.
\item $\mu^*_2(\gamma,\kappa):=\sup \{ \mu_{T}^*(\gamma,\kappa) \; |
\;T\in L_{\kappa^+,\omega}\}$\footnote{Note that similarly to what we did in 
the previous section  with the 
function $\mu_0(\cdot  ,\cdot )$ above, the functions $\mu^*(\gamma,\kappa)$  and
 $\mu^*(\lambda,\kappa)$ are  different objects, we distinguish between them by using
different arguments.}.
\end{enumerate}
\end{definition}

The improvement in comparison  to what we  have seen before is that instead of 
limiting attention
to formulas with the order-property to be from $L_{\lambda^+,\omega}$ we 
consider what may look as a
weaker order-property, by considering formulas with the order property to be 
from the logic
$L_{\infinity ,\omega}$ (with rank bounded by $\gamma$).

\begin{definition}
Let  $T,<,<^P,P$  be as in Definition \ref{delta1}.  
For an \underline{ordinal}  $\gamma>\kappa$ let 
 \begin{enumerate}
\item $\delta_2(\lambda,\gamma,\kappa):=\min\{\delta : \;\;\Gamma$ is a set of $T$-types
,$\;|\Gamma|\leq\kappa, |T|
\leq\lambda\; \\$
\underline{if} $\;\forall\delta'<\delta\;\;\exists M\in EC(T,\Gamma)\;\;$ with
$otp(P^M,<^{P^M})<\gamma$  and  $otp(M-P,<)\geq\delta'\;$, \underline{then} 
$\exists N\in EC(T,\Gamma) $  s.t.  $otp(P^N,<^{P^N})\in On\cap\kappa^+$  and 
$(N-P^N,<^N)$ is not well ordered\}.

\item $\mu_2(\lambda,\gamma,\kappa):=\min\{\mu : \;\;\Gamma$ is a set of $T$-types
,$\;|\Gamma|\leq\kappa, |T|
\leq\lambda\;$  \underline{if} $\exists M\in EC(T,\Gamma)\;\;
 \|M\|\geq\mu $  with $otp(P^M,<^M)<\gamma$ \underline{then}
 for every $\chi\geq\kappa\;\;\exists
N\in EC(T,\Gamma) $  of cardinality at least $\chi$ such that
$otp(P^N,<^N)\in On\cap\kappa^+\}$.

\item  When $\lambda=\kappa$ we may omit $\lambda$.  By the discussion after 
Theorem \ref{ineq} this case is interesting enough. 
\end{enumerate}
\end{definition}

The following is an analog of Proposition \ref{mu1=mu*}:

\begin{proposition}  Let  $\kappa$ and $\mu$ be cardinalities,  and let 
$\gamma$ be an ordinal
$\;\geq\kappa$.  Then $(1)\Rightarrow(2)\Rightarrow(3)$ where
\begin{enumerate}
\item  $\mu\geq\mu_2(\gamma,\kappa)$
\item for every $\psi\in L_{\kappa^+,\omega}$, and for every 
$\phi(x;y)\in L_{\infinity,\omega}$ of quantifier depth $\;<\gamma$ 
\underline {if}  $\psi$  has the $(\phi,\mu)$-order property
\underline {then} there exists $\phi'\in L_{\kappa^+,\omega}$ such that $\psi$ has the
$(\phi',\infinity)$-order property.
\item $\mu\geq\beth_{\gamma}$.
\end{enumerate}
\end{proposition}

\begin{theorem}\label{eq2}
For every $\kappa$  and every ordinal  $\gamma\geq\kappa$   we have  
$\mu_2(\gamma,\kappa)=\beth_{\delta_2(\gamma,\kappa)}$.
\end{theorem}

Theorem \ref{eq2} will be proved in the next section.\\

The following theorem connects  $\delta_2$ to the Galvin--Hajnal rank and 
provides a lower bound for $\delta_2(\gamma,\kappa)$: 

\begin{theorem}\label{lb} {\rm (A lower bound):}  Suppose $\kappa$ is an uncountable regular 
cardinality.
Let $J$ be the ideal of nonstationary subsets of $\kappa$.  For every ordinal $\gamma >
 \kappa$ we have 
$\|\gamma\|_J<\delta_2(\gamma ,\kappa )$, when $\|\gamma \|_J$ is the Galvin--Hajnal 
rank of the 
constant function $f\;:\;\kappa\rightarrow
\gamma+1$ whose value is $\gamma$.
\end{theorem}

Instead of proving the above theorem, we prove a more general result.  
It turns out that the ideal
 $J$ of nonstationary subsets can be replaced by almost any other ideal satisfying rather
weak conditions:

\begin{theorem}\label{lowerbound} {\rm (A better lower bound):}  
Suppose $J$ is an $\aleph_1$-complete
ideal on $\kappa$ such that \\
(*)  $J$ as an ideal is generated by $\leq\kappa$ sets or at least we have\\
(**)  there exists a model ${\cal B}$ (of an expansion of set theory) with universe 
$\kappa , \;\;
|L({\cal B})|\leq\kappa$  and $\psi(P)\in L_{\kappa^+,\omega}$,  when  
$L=L({\cal B})\cup \{P\}$,  $P$ is a unary predicate; having the following property:\\

$\bigotimes_J\;$ for every  $A\ci\kappa$,  we have that  
$ A\in J\iff\langle {\cal B},A\rangle\models\psi(P)$\\ 

or at least \\

$\bigotimes^-_J\;$ for every  $A\ci\kappa$,  we have that  
$ A\in J$ iff for some $A'$ we have that $A\subseteq A'\in J,\;
\langle {\cal B},A'\rangle\models\psi(P)$\\

then
for every ordinal $\gamma>\kappa$  we  have that 
$\|\gamma\|_J < \delta_2(\gamma,\kappa)$.

\end{theorem}

\begin{remark}{\em
\begin{enumerate}
\item  One way to see that Theorem \ref{lb} is a special case of Theorem \ref{lowerbound}
is by using the same argument. Another formal argument (using the statement
of \ref{lowerbound}) we can take 
${\cal B}:=\langle\kappa,<\rangle$ and $\psi(P)$ will say that $\{x\;:\; P(x)\}$ 
is a closed unbounded set.  This satisfy $\bigotimes^-_J\;$ but not $\bigotimes_J$.
\item Note that $\bigotimes^-_J\;$ is equivalent to:  for some  
$\psi(P,{\bar R})\in L_{\kappa^+,\omega}$, we have that 
$ A\in J\iff (\exists {\bar R})\langle{\cal B},A\rangle\models\psi(P,{\bar R})$.

\end{enumerate}

}\end{remark}

{\em Proof.} Let $\gamma^*:=\|\gamma\|_J$, and let  
$ds(\gamma^*)$ stand for the set \\ 
$\{\nu | \nu\mbox { is strictly decreasing sequence of ordinals } <\gamma^*\}$.
There exists a family of functions 
$\{f_{\eta} : \kappa\rightarrow On \;|\; \eta\in ds(\gamma^*) \}$
with the following properties:\begin{enumerate}
\item  $f_{\langle\rangle}$ is constantly  $\gamma$.
\item if  $\eta\^\;i\in ds(\gamma^*)$  then  $f_{ \eta\^{\;}i } <_{J} f_{\eta}$,  and
for every  $\zeta<\kappa$ we have  \\ 
$\neg [f_{\eta\^\;i } (\zeta)  < f_{\eta} (\zeta)] \Rightarrow 
f_{\eta\^\;i}(\zeta )=f_{\eta } (\zeta)=0$.
\item if  $\eta\neq\langle\rangle$  then $\forall \zeta<\kappa [f_{\eta}(\zeta)<\gamma]$.
\item $\eta\^\;i\in ds(\gamma^*)\Rightarrow  \|f_{\eta\^\;i}\|_J\geq i$.
\item $\|f_{\langle\rangle}\|_J=\gamma^*$.

\end{enumerate}
{\sf This is possible:}  Define the function $f_{\eta}$  by induction on  $l(\eta)$:\\
For $l(\eta)=0$; Let $f_{\langle\rangle}$ be the constant function as 
in requirement (1).\\  
For $l(\eta)>0$; If $\eta\^\;i\in ds(\gamma^*)$ then  $f_{\eta}$  is defined, and
 by the inductive hypothesis we have that $\|f_{\eta\^\;i}\|_J>i$   
(as $\eta\^\;i\in ds(\gamma^*)$  and  $\|f_{\langle\rangle}\|_J=\gamma^*$), 
by the definition of the Galvin--Hajnal rank there exists $f'<_Jf_{\eta}$  such that
$\|f'\|_J\geq i$.  Now for $\zeta<\kappa$ let
\[ f_{\eta\^\;i } (\zeta ):=\left\{\begin{array}{ll}
f'(\zeta) & \mbox{if $f'(\zeta)<f_{\eta}(\zeta)$}\\
0 & \mbox{otherwise.} \end{array} \right. \]

{\sf This is enough:} Denote by $\overline{f}$  the sequence 
$\langle f_{\eta} \;:\; \eta\in 
ds(\gamma^*)\rangle$.  \\  Let $\chi^*$  be a sufficiently large regular cardinal such that $H(\chi^*)$
conatins all relevant sets and the structure $\langle H(\chi^*),\in \rangle $
 reflects all relevant
properties.
Let  ${\mathfrak C}:=\langle H(\chi^*),\in,<_{\chi^*}^*,\overline{f},\kappa,D,
{\cal B},\psi(\cdot),J,
Q,P,i\rangle_{i\leq\kappa}$ ,  when  $<_{\chi^*}^*$ is a 
well ordering of the set $ H(\chi^*)$, 
$P$ is the unary predicate $\{i\;:\; i<\gamma \}$,  $Q$ is the unary predicate 
interpreted by
the set $\{j\;:\; j<\gamma^*\}$, $D$ interpreted by $ds(\gamma^*)$.  
Let  $T:=Th({ \mathfrak C})$,
and $\Gamma$  is a set of types consisting  only of the following type -- 
$\{x\in\kappa\wedge
x \neq i\;:\; i<\kappa\}$.  Suppose $N\in EC(T,\Gamma)$ is such
that $(P^N,<^{P^N})$  is well ordered, and we will show that $(Q^N,<)$  is well ordered.

W.l.o.g. we may assume that $A^N:=\{x\in N \;:\; N\models rk(x)\in P \}$ 
is a transitive set and
$\in ^N\restriction A^N=\in\restriction A^N$  (by taking the Mostowski's collapse). 
So $P^N=\gamma'$
for some $\gamma'$,  and since $N$ omits the type in $\Gamma$ we have  
$\kappa^N=\kappa$, since the universe of ${\cal B}$ is $\kappa$ we have 
${\cal B}^N={\cal B}$.  
So necessarily  $\psi(\cdot)^N=\psi(\cdot)$.

By the axioms of $T$ it follows that \\
(*)  if  $\; N\models \eta\^\;x\in D$  then  $N\models (\exists X\subseteq
\{   i<\kappa \;:\; 
f_{\eta}(i)\leq f_{\eta\^\; x}(i)  \})\;
{\cal B}\models\psi(X)$.  Now since $\eta\^\; x\in D^N$,  we have 
$f^N_{\eta},\; f^N_{\eta\^\; x}\in A^N$.  So by the functions from $\kappa$  
into $\gamma '$. Also 
$\psi^N=\psi$, by absoluness we have 
$N\models (\exists X\subseteq
\{   i<\kappa \;:\; 
f_{\eta}(i)\leq f_{\eta\^\; x}(i)  \})\;
{\cal B}\models\psi(X)$.  
So by (**) we have\\
(*)  ${\eta\^\; x}\in D^N\Rightarrow f_{\eta\^\; x} <_J f_{\eta}$.

Now if  $(Q^N,<)$ is not well ordered then we can find 
$\{x_n \in Q^N \;:\; N\models [x_{n+1}<x_n]\}$.
From $T$'s axioms it follows that there 
are  $\{y_n \;:\; n<\omega\}$ such that 
$y_0=\langle \rangle,\;\; y_{n+1}=y_n\^\; x_n\in D$  for all $n<\omega$.  So we have that  
$\{f_{y_n}:\kappa\rightarrow \gamma'\;|\;n<\omega\}$  and  for every $n<\omega$
$f_{y_{n+1}}<_J f_{y_n}$  (in $V$).  Since  $J$ is an
$\aleph_1$-complete ideal we have a contradiction.
We have shown that there exists a pair $T,\Gamma$ of the aproppriate cardinalities such 
that \begin{enumerate}
\item $N\in EC(T,\Gamma), \; (P^N,<)$ is well ordered $\Rightarrow(Q^N,<)$
is well ordered.
\item there is  $N\in EC(T,\Gamma$) with $otp(P^N,<)\leq \gamma$  and  $(Q^N,<)$ 
of order type 
$\gamma^*$  (take $N={\mathfrak C}$).
\end{enumerate}
This establishes that $\gamma^*< \delta_2(\gamma,\kappa)$.
\hfill
$\;\;\;\square_{\ref{lowerbound}} $

\newpage
\section{Concluding Remarks}

It is natural to ask whether the lower bound from Theorem \ref{lowerbound} is equal
to the one in Fact \ref{covfact}.  The following seems to be a reasonable

\begin{conjecture}  For cardinalities $\kappa$ of uncounable cofinality
and $\lambda$ such that $2^{\kappa}<\lambda$  we have 
$\delta_1(\lambda,\kappa)=(cov(\lambda,\kappa)+2^{\kappa})^+$.

\end{conjecture}
\begin{remark}{\rm \begin{enumerate}
\item
Notice that when $\kappa$ is strong limit  singular of cofinality $\aleph_0$
then the conjecture holds.
\item
Why $2^{\kappa}<\lambda$ -- See Barwise-Kunen for independence results.
\item
The conjecture can  
to large extent  be
traslated to a one on $pcf$; it is evident that 
e.g. (***) below is a sufficient condition:

(***) for any set $\mathfrak {a}$ of $\le \kappa$ regular cardinals which are $>2^{\kappa}$
the set $pcf(\mathfrak {a})$ has cardinality at most $\kappa$, or at least
the set  $pcf_{\aleph_1-complete}(\mathfrak {a})$ 
has cardinality at most $\kappa$.

This is because  by \cite{Sh g} II 5.4 if $2^{\kappa}<\lambda$ then $\mu=cov(\lambda,\kappa)$ 
is the first
$\mu$ such that if $\{\lambda_i\;:\; i<\kappa\}$ is a set of
regular cardinalities in the interval $(2^{\kappa},\lambda)$ and $J$
 is an $\aleph_1$-complete ideal on $\kappa$ and $cf(\prod_{i<\kappa}\lambda_i ,<_J)$ 
is well defined then it is $\leq\mu$.

The problem is that the ideal may not satisfy even $\bigotimes^-_J$.  However by
\cite{Sh g} VII, 2.6 the ideal $J$ is generated by a family of 
$\leq|pcf\{\lambda_i\;:\; i<\kappa\}|$  sets
and even by a family of just
$\leq|pcf_{\aleph_1-complete}\{\lambda_i\;:\; i<\kappa\}|$  sets
 , so we have $|pcf(\{\lambda_i\;:\; i<\kappa\}|\leq\kappa
\Rightarrow \bigotimes_J\;$
\end{enumerate} }
\end{remark}


\newpage

\begin{thebibliography}{99}

\bibitem[BaKu]{BaKu} K. J. Barwise and K. Kunen, 
\newblock {Hanf numbers for fragments of $L_{\infinity,\omega }$},
\newblock {\em Israel J. of Mathematics}, {\bf 10}, 306--320, 1971

\bibitem[Ch]{Ch} C. C. Chang,
\newblock {Some remarks on the model theory of infinitary languages},
\newblock {\bf The Syntax and Semantics of Infinitary Languages}, 
Springer-Verlag Lecture Notes
in Mathematics {\bf 72}, edited by J. Barwise (1968)  pages 36--63. 

\bibitem[Gr]{Gr} Rami Grossberg, {\bf Classification theory for non elementary classes}, 
Ph.D. Thesis 
(in Hebrew), The Hebrew University of Jerusalem 1986.

\bibitem[GrSh]{GrSh} Rami Grossberg and Saharon Shelah,
\newblock {On the number of nonisomorphic models of an infinitary theory which
  has the infinitary order property. I},
\newblock {\em {The Journal of Symbolic Logic}}, {\bf 51}:302--322, 1986.
  ---  {\bf MR:}~87j:03037, (03C45)

\bibitem[GrSh 238]{GrSh238}Rami Grossberg and Saharon Shelah,
\newblock {A nonstructure theorem for an infinitary theory which has the
  unsuperstability property},
\newblock {\em {Illinois Journal of Mathematics}}, {\bf 30}:364--390, 1986.
\newblock {Volume dedicated to the memory of W.W.~Boone; ed. Appel, K., Higman,
  G., Robinson, D. and Jockush, C.}
  ---  {\bf MR:}~87j:03036, (03C45)

\bibitem[HrSh 334]{HrSh}Ehud Hrushovski and Saharon Shelah,
\newblock {Stability and omitting types},
\newblock {\em {Israel Journal of Mathematics}}, {\bf 74}:289--321, 1991.
  ---  {\bf MR:}~92m:03049, (03C45)

\bibitem[Ma]{Ma} Menachem Magidor, Chang's conjecture and powers of singular cardinals, 
{\em {The Journal of Symbolic Logic}}, {\bf 42}: 272--276, 1977.

\bibitem[MaSh 285]{MaSh 285} Michael Makkai and Saharon Shelah,
\newblock {Categoricity of theories in $L_ {\kappa\omega},$ with $\kappa$ a
  compact cardinal},
\newblock {\em {Annals of Pure and Applied Logic}}, {\bf 47}:41--97, 1990.
  ---  {\bf MR:}~92a:03054, (03C75)

\bibitem[Mo]{Mo}  Michael Morley, Omitting classes of elements, 
{\bf The Theory of Models}, edited by 
Addison, Henkin and Tarski,  
North-Holland Publ Co  (1965) pages 265--273.


\bibitem[Sh c]{bible}Saharon Shelah.
\newblock {\em {Classification theory and the number of nonisomorphic models}},
  volume~92 of {\em {Studies in Logic and the Foundations of Mathematics}}.
\newblock {North-Holland Publishing Co., Amsterdam, xxxiv+705 pp}, 1990.

\bibitem[Sh g]{Sh g}Saharon Shelah.
\newblock {\em {Cardinal Arithmetic}}, volume~29 of {\em {Oxford Logic
  Guides}}.
\newblock {Oxford University Press}, 1994.


\bibitem[Sh 12]{Sh12}Saharon Shelah.
\newblock {The number of non-isomorphic models of an unstable first-order
  theory}.
\newblock {\em {Israel Journal of Mathematics}}, {\bf 9}:473--487, 1971.
  ---  {\bf MR:}~43-4652.

\bibitem[Sh 16]{Sh16}Saharon Shelah,
\newblock {A combinatorial problem; stability and order for models and theories
  in infinitary languages},
\newblock {\em {Pacific Journal of Mathematics}}, {\bf 41}:247--261, 1972.
  ---  {\bf MR:}~46:7018, (02H10)

\bibitem[Sh 87a]{Sh 87a}Saharon Shelah.
\newblock {Classification theory for nonelementary classes, I. The number of
  uncountable models of $\psi \in L_{\omega _{1},\omega }$. Part A}.
\newblock {\em {Israel Journal of Mathematics}}, {\bf 46}:212--240, 1983.

\bibitem[Sh 87b]{Sh 87b}Saharon Shelah.
\newblock {Classification theory for nonelementary classes, I. The number of
  uncountable models of $\psi \in L_{\omega _{1},\omega }$. Part B}.
\newblock {\em {Israel Journal of Mathematics}}, {\bf 46}:241--273, 1983.

\bibitem[Sh 88]{Sh 88}Saharon Shelah.
\newblock {Classification of nonelementary classes. II. Abstract elementary
  classes}.
\newblock In {\em {Classification theory (Chicago, IL, 1985)}}, volume 1292 of
  {\em {Lecture Notes in Mathematics}}, pages 419--497. {Springer, Berlin},
  1987.
\newblock {Proceedings of the USA--Israel Conference on Classification Theory,
  Chicago, December 1985; ed. Baldwin, J.T.}

\bibitem[Sh 300]{Sh 300}Saharon Shelah.
\newblock {Universal classes}.
\newblock In {\em {Classification theory (Chicago, IL, 1985)}}, volume 1292 of
  {\em {Lecture Notes in Mathematics}}, pages 264--418. {Springer, Berlin},
  1987.
\newblock {Proceedings of the USA--Israel Conference on Classification Theory,
  Chicago, December 1985; ed. Baldwin, J.T.}

\bibitem[Sh 299]{Sh 299}Saharon Shelah.
\newblock {Taxonomy of universal and other classes}.
\newblock In {\em {Proceedings of the International Congress of Mathematicians
  (Berkeley, Calif., 1986)}}, volume~1, pages 154--162. {Amer. Math. Soc.,
  Providence, RI}, 1987.
\newblock {ed. Gleason, A.M.}

\bibitem[Sh 78]{Sh78}Saharon Shelah.
\newblock {Hanf number of omitting type for simple first-order theories}.
\newblock {\em {The Journal of Symbolic Logic}}, {\bf 44}:319--324, 1979.


\bibitem[Sh 48]{Sh 48}Saharon Shelah.
\newblock {Categoricity in $\aleph _{1}$ of sentences in $L_{\omega _{1},\omega
  }(Q)$}.
\newblock {\em {Israel Journal of Mathematics}}, {\bf 20}:127--148, 1975.

\bibitem[Sh 200]{Sh 200}Saharon Shelah.
\newblock {Classification of first order theories which have a structure
  theorem}.
\newblock {\em {American Mathematical Society. Bulletin. New Series}}, {\bf
  12}:227--232, 1985.

\end{thebibliography}
\end{document}

--ELM906474423-15724-0_--